\documentclass[12pt, dvipdfmx]{article}
\usepackage[mathscr]{eucal}
\usepackage{amssymb}
\usepackage{latexsym}
\usepackage{amsthm}
\usepackage{amsmath}
\usepackage[dvips]{graphicx}
\usepackage{psfrag}
\usepackage{a4wide}
\usepackage{tikz-cd}
\usepackage{amscd}
\usepackage{url}

\theoremstyle{plain}
\newtheorem{theorem}{Theorem}[section]
\newtheorem{proposition}[theorem]{Proposition}
\newtheorem{lemma}[theorem]{Lemma}
\newtheorem{corollary}[theorem]{Corollary}

\theoremstyle{definition}
\newtheorem*{definition}{Definition}

\theoremstyle{remark}

\newcommand{\C}{{\mathbb C}}
\newcommand{\D}{{\mathbb D}}

\newcommand{\cS}{{\mathcal S}}

\newcommand{\cH}{{\mathcal H}}

\newcommand{\la}{\langle}
\newcommand{\ra}{\rangle}

\newcommand{\lam}{\lambda}

\begin{document}

\title{
Operator theory induced by powers of the de Branges-Rovnyak kernel 
and its application\footnote{
This paper has been accepted by Canadian Mathematical Bulletin, 
in which the new title is ``Exponentials of de Branges-Rovnyak kernels''. }}
\author{
{\sc Shuhei KUWAHARA}\\
[1ex]
{\small Sapporo Seishu High School, Sapporo 064-0916, Japan}\\
{\small 
{\it E-mail address}: {\tt s.kuwahara@sapporoseishu.ed.jp}}\\
and\\ 
{\sc Michio SETO}\\
[1ex]
{\small National Defense Academy,  
Yokosuka 239-8686, Japan} \\
{\small 
{\it E-mail address}: {\tt mseto@nda.ac.jp}}
}

\date{}

\maketitle

\begin{abstract}
In this note, 
we give a new property of de Branges-Rovnyak kernels.      
As the main theorem, it is shown that 
the exponential of de Branges-Rovnyak kernel is strictly positive definite
if the inner part of the corresponding Schur class function is nontrivial.
\end{abstract}

\begin{center}
keywords: de Branges-Rovnyak kernel, positive definite kernel
\end{center}

\section{Introduction}
Let $\D$ be the open unit disk in the complex plane $\C$, and
let $H^{\infty}$ be the Banach algebra consisting of all bounded analytic functions on $\D$.
Then, we set
\[
\cS=\{\varphi\in H^{\infty}:|\varphi(\lam)|\leq 1\ (\lam \in \D)\}, 
\]
and which is called the Schur class. 
For any function $\varphi$ in $H^{\infty}$, 
it is well known that $\varphi$ belongs to $\cS$ if and only if 
\begin{equation}\label{eq:1-1}
\dfrac{1-\overline{\varphi(\lam)}\varphi(z)}{1-\overline{\lam}z}
\end{equation}
is positive semi-definite.
This equivalence relation based on the properties of the Szeg\"{o} kernel 
is crucial in the operator theory on the Hardy space over $\D$, in particular, 
theories of Pick interpolation, de Branges-Rovnyak spaces and sub-Hardy Hilbert spaces 
 (see Agler-McCarthy~\cite{AM}, Ball-Bolotnikov~\cite{BB}, Fricain-Mashreghi~\cite{FM} and Sarason~\cite{Sarason}). 
The kernel (\ref{eq:1-1}) is called the de Branges-Rovnyak kernel. 

Before introducing our study, 
we should mention that not only the original de Branges-Rovnyak kernel but
also its variants have been studied by a number of authors. 
For example, Zhu~\cite{Zhu1,Zhu2} initiated the study on 
the kernel 
\begin{equation}\label{eq:1-2}
\dfrac{1-\overline{\varphi(\lam)}\varphi(z)}{(1-\overline{\lam}z)^2}
\end{equation}
in the Bergman space over $\D$. 
The reproducing kernel Hilbert space induced by the kernel 
(\ref{eq:1-2}) is called a sub-Bergman Hilbert space 
(see also Abkar-Jafarzadeh~\cite{AJ}, Ball-Bolotnikov~\cite{BB1}, Chu~\cite{Chu}, 
Nowak-Rososzczuk~\cite{NR} and Sultanic~\cite{Sultanic}). 
Further, powers of the de Branges-Rovnyak kernel
\begin{equation}\label{eq:1-3}
\left(\dfrac{1-\overline{\varphi(\lam)}\varphi(z)}{1-\overline{\lam}z}\right)^n\quad (n\in {\mathbb N})
\end{equation}
are naturally obtained from the theory of hereditary functional calculus for weighted Bergman spaces on $\D$ (see Example 14.48 in \cite{AM} for the case where $n=2$) 
and have appeared also in p. 3672 of Jury~\cite{J}.  

Now, the purpose of this paper is to study the structure of the kernel 
\begin{equation}\label{eq:1-4}
\exp\left(t\dfrac{1-\overline{\varphi(\lam)}\varphi(z)}{1-\overline{\lam}z}\right)\quad (t>0).
\end{equation}
Note that our kernel (\ref{eq:1-4}) is 
obtained by binding all kernels in (\ref{eq:1-3}) together. 
Thus, we expect that new properties of the de Branges-Rovnyak kernel (\ref{eq:1-1}) are drawn out 
from our kernel (\ref{eq:1-4}).  
In fact, as the main result, we will show 
that the exponential of the de Branges-Rovnyak kernel is strictly positive definite
if the inner part of $\varphi$ is nontrivial. 

Here, we shall give some remarks on strictly positive definite kernels.  
In general, it is not difficult to construct positive semi-definite kernels. 
On the other hand, for strictly positive definite kernels, nontrivial methods depending on each case are often needed 
(for example, see Micchelli~\cite{M}).    
Moreover, it might be worth while mentioning that strictly positive definite kernels have received attention 
in machine learning (see Rasmussen-Williams~\cite{RW}). 

This paper is organized as follows. 
In Section 2, basic properties of 
the reproducing kernel Hilbert space $\exp \cH_t(\varphi)$ constructed from our kernel (\ref{eq:1-4}) are given. 
In Section 3, unbounded multipliers on $\exp \cH_t(\varphi)$ are introduced and studied. 
In Section 4, main results are given.  

The problem discussed in this paper was obtained in conversation about machine learning with 
Professor Kohtaro Watanabe (National Defense Academy). 
The authors would like to express gratitude to him. 
This research was supported by JSPS KAKENHI Grant Number 20K03646.

\section{Preliminaries}
For $t>0$, 
let $\cH_t(\varphi)$ denote the reproducing kernel Hilbert space 
with kernel
\[
tk^{\varphi}(z,\lam)
=t\dfrac{1-\overline{\varphi(\lam)}\varphi(z)}{1-\overline{\lam}z}\quad (\varphi \in \cS), 
\] 
and we will use notations $tk_{\lam}^{\varphi}(z)=tk^{\varphi}(z,\lam)$ and $\cH(\varphi)=\cH_1(\varphi)$. 
Then, since
\[
\la tk_{\lam}^{\varphi},tk_z^{\varphi} \ra_{\cH_t(\varphi)}
=tk^{\varphi}(z,\lam)=t^{-1}\la tk_{\lam}^{\varphi},tk_z^{\varphi} \ra_{\cH(\varphi)}, 
\]
the trivial linear mapping $f\mapsto f$ 
from $\cH(\varphi)$ onto $\cH_t(\varphi)$ is bounded and invertible. 
Particularly, $\cH_t(\varphi)=\cH(\varphi)$ as vector spaces. 
In this section, we construct the exponential of $\cH_t(\varphi)$ and give its basic properties. 
The contents of this section are well known to specialists. 
For example, see Exercise (k) in p. 320 of Nikolski~\cite{Nik} 
and Chapter 7 in Paulsen-Raghupathi~\cite{PR}. 
However, we give the details for the sake of readers. 

\subsection{Construction of \mathversion{bold}$\exp \cH_t(\varphi)$}
Let $\cH_t(\varphi)^n$ be the reproducing kernel Hilbert space obtained 
by the pull-back construction with the $n$-fold tensor product space 
\[
\cH_t(\varphi)^{\otimes n}=\cH_t(\varphi)\otimes \cdots \otimes \cH_t(\varphi)
\]
and the $n$-dimensional diagonal map 
\[
\Delta_n:\D\to \D^n,\ \lam\to (\lam,\ldots,\lam)
\]
(for the pull-back construction, see Theorem 5.7 in \cite{PR}). 
We note that $
(tk_{\lam}^{\varphi})^{\otimes n}\circ \Delta_n=(tk_{\lam}^{\varphi})^n$ 
is the reproducing kernel of $\cH_t(\varphi)^n$.  
Let $\oplus_{n=0}^{\infty}\cH_t(\varphi)^n$ denote the Hilbert space with the inner product
\[
\la (f_0,f_1,\ldots)^{\top},(g_0,g_1,\ldots)^{\top} \ra_{\oplus_{n=0}^{\infty}\cH_t(\varphi)^n}
=\sum_{n=0}^{\infty}\dfrac{1}{n!}\la f_n,g_n \ra_{\cH_t(\varphi)^n},
\]
where we set $\cH_t(\varphi)^0=\C$. 
Moreover, we define linear map $\Gamma$ as follows: 
\[
\Gamma: 
\begin{pmatrix}
f_0\\
f_1\\
\vdots
\end{pmatrix}
\mapsto \sum_{n=0}^{\infty}\dfrac{1}{n!}f_n\quad 
\left(
\begin{pmatrix}
f_0\\
f_1\\
\vdots
\end{pmatrix}
\in\oplus_{n=0}^{\infty}\cH_t(\varphi)^n
\right).
\] 

\begin{proposition}\label{prop:3-1}
The following statements hold:
\begin{enumerate}
\item $\Gamma$ is a map from $\oplus_{n=0}^{\infty}\cH_t(\varphi)^n$ 
to $\mbox{\rm Hol}(\D)$. 
\item $\ker \Gamma$ is closed. 
\end{enumerate}
\end{proposition}

\begin{proof}
For any $F=(f_0,f_1,\ldots)^{\top}$ in $\oplus_{n=0}^{\infty}\cH_t(\varphi)^n$, 
we have 
\begin{align*}
\left| \sum_{\ell=n+1}^m\dfrac{1}{\ell !}f_{\ell}(\lam) \right|
&\leq 
\sum_{\ell=n+1}^m\left|\dfrac{1}{\ell !}f_{\ell}(\lam)\right|\\
&\leq \sum_{\ell=n+1}^m\dfrac{1}{\ell !}\|f_{\ell}\|_{\cH_t(\varphi)^{\ell}}\|(tk_{\lam}^{\varphi})^{\ell}\|_{\cH_t(\varphi)^{\ell}}\\
&\leq \left(\sum_{\ell =n+1}^m\dfrac{1}{\ell !}\|f_{\ell}\|_{\cH_t(\varphi)^{\ell}}^2\right)^{1/2}
\left(\sum_{\ell=n+1}^m\dfrac{1}{\ell !}\|(tk_{\lam}^{\varphi})^{\ell}\|_{\cH_t(\varphi)^{\ell}}^2\right)^{1/2}\\
&= \left(\sum_{\ell=n+1}^m\dfrac{1}{\ell !}\|f_k\|_{\cH_t(\varphi)^{\ell}}^2\right)^{1/2}
\left(\sum_{\ell=n+1}^m\dfrac{1}{\ell !}\|tk_{\lam}^{\varphi}\|_{\cH_t(\varphi)}^{2\ell}\right)^{1/2}.
\end{align*}
Hence, 
\[
\sum_{n=0}^{\infty}\dfrac{1}{n!}f_n(\lam)
\]
converges uniformly on any compact subset of $\D$. This concludes (1). 
Next, suppose that $F_{\ell}=(f_0^{(\ell)},f_1^{(\ell)},\ldots)^{\top}$ belongs to $\ker \Gamma$ and 
$F_{\ell}$ converges to $F$ in $\oplus_{n=0}^{\infty}\cH_t(\varphi)^n$. 
Then, for sufficiently large $L$, we have 
\[
\sum_{n=0}^{\infty}\dfrac{1}{n!}\|f_n^{(\ell)}-f_n\|_{\cH_t(\varphi)^n}^2
=\|F_{\ell}-F\|_{\oplus_{n=0}^{\infty}\cH_t(\varphi)^n}^2
<1\quad (\ell \geq L).
\]
Hence we have
$\|f_n^{(\ell)}-f_n\|_{\cH_t(\varphi)^n}<\sqrt{n!}$ for every $n\geq 0$ if $\ell\geq L$. 
It follows from this inequality that, for any $\ell \geq L$,  
\begin{align*}
\left|\dfrac{1}{n!}f_n^{(\ell)}(\lam)\right|
&\leq \dfrac{1}{n!}\|f_n^{(\ell)}\|_{\cH_t(\varphi)^n}\|(tk_{\lam}^{\varphi})^n\|_{\cH_t(\varphi)^n}\\
&\leq \dfrac{(\|f_n\|_{\cH_t(\varphi)^n}+\sqrt{n!})\|tk_{\lam}^{\varphi}\|_{\cH_t(\varphi)}^n}{n!}\\
&\leq \dfrac{(M+1)\|tk_{\lam}^{\varphi}\|_{\cH_t(\varphi)}^n}{\sqrt{n!}},
\end{align*}
where we set $M=\sup_{n\geq 0} (\|f_n\|_{\cH_t(\varphi)}^2/n!)^{1/2}$. 
Then, by the ratio test, 
\[
\sum_{n=0}^{\infty}\dfrac{(M+1)\|tk_{\lam}^{\varphi}\|_{\cH_t(\varphi)}^n}{\sqrt{n!}}
\]
is finite. 
Hence, by the Lebesgue dominated convergence theorem and the assumption that 
$F_{\ell}=(f_0^{(\ell)},f_1^{(\ell)},\ldots)^{\top}$ belongs to $\ker \Gamma$, we have
\[
\sum_{n=0}^{\infty}\dfrac{1}{n!}f_n(\lam)
=\sum_{n=0}^{\infty}\lim_{\ell \to \infty}\dfrac{1}{n!}f_n^{(\ell)}(\lam)
=\lim_{\ell \to \infty}\sum_{n=0}^{\infty}\dfrac{1}{n!}f_n^{(\ell)}(\lam)=0.
\]
This concludes (2). 
\end{proof}

By Proposition \ref{prop:3-1}, the pull-back construction can be applied to $\Gamma$. 
\begin{definition}
We define $\exp {\cH_t(\varphi)}$ as the reproducing kernel Hilbert space obtained by the pull-back construction 
with the linear map 
\[
\Gamma: \oplus_{n=0}^{\infty}\cH_t(\varphi)^n\to \mbox{\rm Hol}(\D).
\] 
\end{definition}

\subsection{Basic properties of \mathversion{bold}$\exp \cH_t(\varphi)$}
We summarize basic properties of $\exp {\cH_t(\varphi)}$.  

\begin{proposition}
$\exp {\cH_t(\varphi)}$ is a reproducing kernel Hilbert space consisting of holomorphic functions on $\D$. 
More precisely, for any $f$ in $\exp\cH_t(\varphi)$, 
there exists a vector $(f_0,f_1,\ldots,)^{\top}$ in $ \oplus_{n=0}^{\infty}\cH_t(\varphi)^n$ 
such that 
\[
f=\sum_{n=0}^{\infty}\dfrac{1}{n!}f_n
\]
converges uniformly on any compact subset of $\D$. 
Moreover, 
\begin{enumerate}
\item the following norm estimate holds:
\[
\|f\|_{\exp\cH_t(\varphi)}^2\leq \sum_{n=0}^{\infty}\dfrac{1}{n!}\|f_n\|_{\cH_t(\varphi)^n}^2,
\] 
\item the reproducing kernel of $\exp \cH_t(\varphi)$ is 
\[
\sum_{n=0}^{\infty}\dfrac{1}{n!}(tk_{\lam}^{\varphi})^n= \exp tk_{\lam}^{\varphi},
\]
that is, 
\[
f(\lam)=\la f, \exp tk_{\lam}^{\varphi}\ra_{\exp\cH_t(\varphi)}
\]
for any $\lam$ in $\D$,
\item the following growth condition holds: 
\[
|f(\lam)|^2 \leq \|f\|_{\exp\cH_t(\varphi)}^2\exp\left(t \dfrac{1-|\varphi(\lam)|^2}{1-|\lam|^2}\right)
\]
for any $\lam$ in $\D$. 
\end{enumerate}
\end{proposition}

\begin{proof}
By the definition of the norm and the inner product of $\exp\cH_t(\varphi)$, we have (1) and (2). 
We shall show (3). By (2) and the Cauchy-Schwarz inequality, we have
\begin{align*}
|f(\lam)|^2
&=|\la f, \exp tk_{\lam}^{\varphi}\ra_{\exp\cH_t(\varphi)}|^2\\
&\leq \|f\|_{\exp\cH_t(\varphi)}^2\cdot \|\exp tk_{\lam}^{\varphi}\|_{\exp\cH_t(\varphi)}^2\\
&=\|f\|_{\exp\cH_t(\varphi)}^2\cdot\exp tk_{\lam}^{\varphi}(\lam)\\
&=\|f\|_{\exp\cH_t(\varphi)}^2\exp\left( t\dfrac{1-|\varphi(\lam)|^2}{1-|\lam|^2}\right).
\end{align*}
Thus we have (3).
\end{proof}

\section{Unbounded multipliers}
We shall investigate into unbounded multipliers of $\exp\cH_t(\varphi)$. 

\begin{lemma}\label{lem:5-3}
Let $\psi$ be a function in $\cH_t(\varphi)$. Then, 
for any function $f$ in $\cH_t(\varphi)^n$, 
$\psi f$ belongs to $\cH_t(\varphi)^{n+1}$. 
\end{lemma}

\begin{proof}
We define bounded linear operator $\tau_{\psi}$ as follows:  
\[
\tau_{\psi}:
\cH_t(\varphi)^{\otimes n} \to\cH_t(\varphi)^{\otimes n+1}, 
\quad F\mapsto \psi \otimes F.
\]
Then, the following diagram commutes:
\[
\begin{CD}
\cH_t(\varphi)^{\otimes n} @>{\tau_{\psi}}>> \cH_t(\varphi)^{\otimes n+1} \\
@V{\Delta_n}VV    @VV{\Delta_{n+1}}V \\
\cH_t(\varphi)^n   @>>{M_{\psi}|_{\cH_t(\varphi)^n}}>  \cH_t(\varphi)^{n+1},
\end{CD}
\]
where $M_{\psi}$ denotes the multiplication operator with symbol $\psi$. 
This concludes the proof. 
\end{proof}

\begin{theorem}
Let $\psi$ be a function in $\cH_t(\varphi)$. Then, 
the multiplication operator $M_{\psi}$ is a densely defined closable linear operator in $\exp\cH_t(\varphi)$. 
\end{theorem}

\begin{proof}
Let $F=(f_0,f_1,\ldots,f_N,0\ldots)^{\top}$ be a vector having finite support in $\oplus_{n=0}^{\infty}\cH_t(\varphi)^n$. 
We set $\Gamma F=f$. 
Then, 
\[\psi f=\psi \sum_{n=0}^N\dfrac{1}{n!}f_n=\sum_{n=0}^N\dfrac{1}{n!}\psi f_n=
\sum_{n=0}^N\dfrac{1}{(n+1)!}(n+1)\psi f_n=\sum_{n=1}^{N+1}\dfrac{1}{n!}n\psi f_{n-1},
\]
where we note that $n\psi f_{n-1}$ belongs to $\cH_t(\varphi)^n$ by Lemma \ref{lem:5-3}. 
Hence, setting 
\[
G=(0, \psi f_0,2\psi f_1,\ldots, (N+1)\psi f_N,0,\ldots)^{\top},
\] 
$G$ belongs to $\oplus_{n=0}^{\infty}\cH_t(\varphi)^n$ and $\Gamma G=\psi f$, that is, $\psi f$ belongs to $\exp\cH_t(\varphi)$. 
Therefore, $M_{\psi}$ is a densely defined linear operator in $\exp\cH_t(\varphi)$. 
Moreover, 
it is easy to see that $M_{\psi}$ is closable. 
\end{proof}

\begin{corollary}\label{cor:4-5}
Let $\psi$ be a function in $\cH_t(\varphi)$. 
Then the adjoint operator $M_{\psi}^{\ast}$ of $M_{\psi}$ is a densely defined closed linear operator in $\exp \cH_t(\varphi)$, and 
every $\exp tk_{\lam}^{\varphi}$ is an eigenfunction of $M_{\psi}^{\ast}$. 
More precisely, 
\[
M_{\psi}^{\ast}\exp tk_{\lam}^{\varphi}=\overline{\psi(\lam)}
\exp tk_{\lam}^{\varphi}.
\]
\end{corollary}

\section{Main results}

Let $X$ be a set. 
A function $k$ on $X\times X$ is called a strictly positive definite kernel on $X$ if 
$k(x,y)=\overline{k(y,x)}$ for any $x$ and $y$ in $X$ and 
\[
\sum_{i,j=1}^nc_i\overline{c_j}k(x_j,x_i)>0
\]
for any $n$ in $\mathbb N$, any $(c_1,\ldots, c_n)^{\top}$ in ${\mathbb C}^n\setminus \{ \mathbf{0} \}$ 
and any $n$ distinct points $x_1,\ldots,x_n$ in $X$. 
For example, it is well known that 
\[
k(z,\lam)=\exp(\overline{\lam}z)
\]
is a strictly positive definite kernel on $\C$. 
In fact, this is the reproducing kernel of the Segal-Bargmann space. 
Now, we note that if $\varphi=z^2$ then 
\[
e^{-1}\exp\left(\dfrac{1-\overline{\varphi(\lam)}\varphi(z)}{1-\overline{\lam}z}\right)
=\exp(\overline{\lam}z).
\]
Motivated by this observation, 
we shall give new examples of strictly positive definite kernels.

Let $\varphi$ be a function in $\cS$. 
For the canonical factorization $\varphi=\alpha z^N BSF$ of $\varphi$ 
(see Section 5 in Chapter II in Garnett~\cite{Garnett}), 
where $\alpha$ is a unimodular constant, 
$B$ is a Blaschke product consisting of nonzero zero points, 
$S$ is a singular inner function and $F$ is an outer function, 
we consider three conditions (C1) $N\geq 2$, 
(C2) $B$ is nontrivial and (C3) $S$ is nontrivial. 
We need the following lemma.

\begin{lemma}\label{lem:6-2}
Let 
$\lam_1,\ldots,\lam_n$ be $n$ distinct points in $\D$. 
Suppose one of {\rm (C1)}, {\rm (C2)} and {\rm (C3)}. 
Then there exists a function $\psi$ in $\cH_t(\varphi)$ such that 
$\psi(\lam_i)\neq \psi(\lam_j)$ $(i\neq j)$.
\end{lemma}

\begin{proof}
Since $\cH_t(\varphi)=\cH(\varphi)$ as vector spaces, 
it suffices to show the statement for $\cH(\varphi)$. 
First, we assume (C1). 
Then, since $\varphi/z$ is in $\cS$ by the Schwarz lemma and $(\varphi/z)(0)=0$, we have
\[
(I-T_{\varphi}T_{\varphi}^{\ast})z=z-T_{\varphi}T_{\varphi/z}^{\ast}T_z^{\ast}z
=z-T_{\varphi}T_{\varphi/z}^{\ast}1=z,
\]
where $T_{\varphi}$ denotes the Toeplitz operator with symbol $\varphi$ on the Hardy space $H^2$ over $\D$. 
Hence $z$ belongs to $\cH(\varphi)$, and we may take $\psi=z$. 

Secondly, we assume (C2). Let $\mu$ be a nonzero zero point of $\varphi$. 
Then, we have  
\[
k_{\mu}^{\varphi}=(I-T_{\varphi}T_{\varphi}^{\ast})k_{\mu}=k_{\mu}. 
\]
Hence $(1-\overline{\mu}z)^{-1}$ belongs to $\cH(\varphi)$, 
and we may take $\psi=(1-\overline{\mu}z)^{-1}$. 

Thirdly, we assume (C3). 
Then, observe that there exists a point $\theta \in [0,2\pi)$ such that
$\varphi(re^{\sqrt{-1}\theta})\to 0$ as $r\uparrow 1$ (see Theorem 6.2 in \cite{Garnett}). 
Without loss of generality, we may assume that $\theta=0$.  
Setting  
\[
\delta_{ij}(r)=\left| k_r^{\varphi}(\lam_i)-k_r^{\varphi}(\lam_j) \right|
\quad \text{and}\quad \delta(r)=\min_{1\leq i<j\leq n} \delta_{ij}(r)\quad (0<r<1),
\]
it suffices to show that $\delta(r)>0$ for any $r$ sufficiently close to $1$.  
For any distinct $\lam_i$ and $\lam_j$, 
we suppose that, 
for any $m$ in $\mathbb N$, there exists 
a real number $r_m \in (1-m^{-1},1)$ such that 
$\delta_{ij}(r_m)=0$. Then, 
we have 
\[
\left|
\dfrac{1}{1-\lam_i}-\dfrac{1}{1-\lam_j}
\right|
=\lim_{m\to \infty}\left|
\dfrac{1-\overline{\varphi(r_m)}\varphi(\lam_i)}{1-r_m\lam_i}
-\dfrac{1-\overline{\varphi(r_m)}\varphi(\lam_j)}{1-r_m\lam_j}
\right|
=\lim_{m\to \infty}\delta_{ij}(r_m)=0.
\]
This concludes that $\lam_i=\lam_j$, however, 
which contradicts the assumption $\lam_i\neq \lam_j$. 
Hence there exists $r_{ij}\in (0,1)$ such that $\delta_{ij}(r)>0$ for any $r\in (r_{ij},1)$. 
Setting $R=\max_{1\leq i<j \leq n} r_{ij}$,  we have $\delta(r)>0$ for any $r\in (R,1)$. Then, 
we may take $\psi=k_r^{\varphi}$. 
\end{proof}

\begin{theorem}\label{thm:6-3}
Let $\varphi$ be a function in $\cS$. 
If $\varphi$ satisfies one of {\rm (C1)}, {\rm (C2)} and {\rm (C3)}, 
then the kernel
\[
k_t(z,\lam)=\exp\left(t\dfrac{1-\overline{\varphi(\lam)}\varphi(z)}{1-\overline{\lam}z}\right)
\quad (t>0)
\]
is strictly positive definite. 
\end{theorem}

\begin{proof}
It suffices to show that 
$\{\exp tk_{\lam_j}^{\varphi}\}_{j=1}^n$ is linearly independent 
for any $n$ in $\mathbb N$ and any $n$ distinct points $\lam_1,\ldots, \lam_n$ in $\D$. 
Suppose that 
\[
\sum_{j=1}^nc_j\exp tk_{\lam_j}^{\varphi}=0
\]
for some $n$ in $\mathbb N$, some $n$ distinct points $\lam_1,\ldots, \lam_n$ in $\D$, 
and some $c_1,\ldots,c_n$ in $\C$. 
Then, for any function $\psi$ in $\cH_t(\varphi)$, 
by Corollary \ref{cor:4-5} and the assumption, we have 
\begin{align*}
\begin{pmatrix}
1 & \cdots & 1\\
\overline{\psi(\lam_1)} & \cdots & \overline{\psi(\lam_n)}\\
\vdots & \vdots & \vdots \\
\overline{\psi(\lam_1)}^{n-1} & \cdots & \overline{\psi(\lam_n)}^{n-1}
\end{pmatrix}
\begin{pmatrix}
c_1\exp tk_{\lam_1}^{\varphi}\\
c_2\exp tk_{\lam_2}^{\varphi}\\
\vdots\\
c_n\exp tk_{\lam_n}^{\varphi}
\end{pmatrix}
&=
\begin{pmatrix}
\sum_{j=1}^nc_j\exp tk_{\lam_j}^{\varphi}\\
\sum_{j=1}^n\overline{\psi(\lam_j)}c_j\exp tk_{\lam_j}^{\varphi}\\
\vdots\\
\sum_{j=1}^nc_j\overline{\psi(\lam_j)}^{n-1}\exp tk_{\lam_j}^{\varphi}
\end{pmatrix}\\
&=
\begin{pmatrix}
\sum_{j=1}^nc_j\exp tk_{\lam_j}^{\varphi}\\
M_{\psi}^{\ast}\sum_{j=1}^nc_j\exp tk_{\lam_j}^{\varphi}\\
\vdots\\
(M_{\psi}^{\ast})^{n-1}\sum_{j=1}^nc_j\exp tk_{\lam_j}^{\varphi}
\end{pmatrix}\\
&=\mathbf{0}.
\end{align*}
Further, by Lemma \ref{lem:6-2}, there exists a function $\psi$ in $\cH_t(\varphi)$ such that 
\[
\prod_{1\leq i<j\leq n}(\psi(\lam_i)-\psi(\lam_j))\neq 0.
\]
Then, the Vandermonde matrix
\[
\begin{pmatrix}
1 & \cdots & 1\\
\overline{\psi(\lam_1)} & \cdots & \overline{\psi(\lam_n)}\\
\vdots & \vdots & \vdots \\
\overline{\psi(\lam_1)}^{n-1} & \cdots & \overline{\psi(\lam_n)}^{n-1}
\end{pmatrix}
\]
is nonsingular. 
Therefore, we have that 
\[
\begin{pmatrix}
c_1\exp tk_{\lam_1}^{\varphi}\\
c_2\exp tk_{\lam_2}^{\varphi}\\
\vdots\\
c_n\exp tk_{\lam_n}^{\varphi}
\end{pmatrix}=\mathbf{0}.
\]
This concludes that $c_1=\cdots=c_n=0$. 
\end{proof}

The well-known fact mentioned at the beginning of 
this section is included in Theorem \ref{thm:6-3}. 
\begin{corollary}\label{cor:4-4}
The kernel function
\[
k(z,\lam)=\exp(\overline{\lam}z)
\]
is strictly positive definite on $\C$. 
\end{corollary}

\begin{proof}
For any $n$ distinct points $\lam_1,\ldots \lam_n$ in $\C$, 
we set $R=\max_{1\leq j \leq n} |\lam_j|+1$. 
Then $\lam_1/R,\ldots \lam_n/R$ are in $\D$.  
Hence, by Theorem \ref{thm:6-3} in the case where $\varphi=z^2$ and $t=R^2$, 
we have 
\begin{align*}
\sum_{i,j=1}^nc_i\overline{c_j}
\exp(\overline{\lam_i}\lam_j)
&=e^{-R^2}\sum_{i,j=1}^nc_i\overline{c_j}
\exp(R^2+\overline{\lam_i}\lam_j)\\
&=e^{-R^2}\sum_{i,j=1}^nc_i\overline{c_j}
\exp(R^2(1+\overline{(\lam_i/R)}(\lam_j/R)))\\
&=e^{-R^2}\sum_{i,j=1}^nc_i\overline{c_j}
\exp\left(R^2\dfrac{1-\overline{\varphi(\lam_i/R)}\varphi(\lam_j/R)}{1-\overline{(\lam_i/R)}(\lam_j/R)}\right)
 >0
\end{align*}
for any $(c_1,\ldots, c_n)^{\top}$ in ${\mathbb C}^n\setminus \{ \mathbf{0} \}$. 
\end{proof}

Although the next result is just a simple consequence of Theorem \ref{thm:6-3}, 
from the viewpoint of the theory of model spaces (see Garcia-Mashreghi-Ross~\cite{GMR}),  
it will be worth while mentioning it as a theorem.  
\begin{theorem}
Let $\varphi$ be an inner function. 
If $\varphi$ is neither a constant nor $e^{i\theta} z$, 
then the kernel
\[
k_t(z,\lam)=\exp\left(t\dfrac{1-\overline{\varphi(\lam)}\varphi(z)}{1-\overline{\lam}z}\right)
\quad (t>0)
\]
is strictly positive definite. 
\end{theorem}

Further, with help of the theory of sub-Hardy Hilbert spaces, we have
\begin{theorem}\label{thm:6-4}
Let $\varphi$ be a function in $\cS$. 
If $\varphi$ is a nonextreme point of the closed unit ball in $H^{\infty}$, 
then the kernel
\[
k_t(z,\lam)=\exp\left(t\dfrac{1-\overline{\varphi(\lam)}\varphi(z)}{1-\overline{\lam}z}\right)
\quad(t>0)
\]
is strictly positive definite. 
\end{theorem}

\begin{proof}
By Theorem (IV-3) in \cite{Sarason} (or Theorem 23.13 in \cite{FM}), the polynomials are dense in $\cH(\varphi)$. 
Hence $z$ belongs to $\cH_t(\varphi)$. 
Then, setting $\psi=z$, the proof of Theorem \ref{thm:6-3} applies to this case. 
\end{proof}

\end{document}